\documentclass[10pt]{amsart}

\usepackage[leqno]{amsmath}
\usepackage{amsthm}
\usepackage{amsfonts}

\theoremstyle{plain}
\newtheorem{lem}{Lemma}[section]
\newtheorem{thm}{Theorem}

\theoremstyle{remark}
\newtheorem{rem}[lem]{Remark}

\newcommand{\ud}{\mathrm{d}}
\newcommand{\Ud}{\mathrm{D}}

\providecommand{\abs}[1]{\left\lvert#1\right\rvert}
\providecommand{\norm}[1]{\left\lVert#1\right\rVert}

\numberwithin{equation}{section}

\DeclareMathOperator{\supp}{supp}
\DeclareMathOperator{\esssup}{ess\, sup}
\DeclareMathOperator{\Div}{div}
\DeclareMathOperator{\Rot}{rot}
\DeclareMathOperator{\Dt}{\frac{\ud}{\ud t}}

\usepackage[T1]{fontenc}
\usepackage[utf8]{inputenc}

\usepackage[margin={3cm,4cm}]{geometry}

\usepackage{caption}
\usepackage{graphicx}
\usepackage{subfig}

\begin{document}

\title[Regularity criteria for NSE involving the pressure]{Regularity criteria of weak solutions to NSE in some bounded domains involving the pressure}

\author{Adam Kubica}
\author{Bernard Nowakowski}
\author{Wojciech M. Zajączkowski}
\thanks{All authors were financially supported by the National Science Centre, under project number NN 201 396937}
\address{Bernard Nowakowski\\ Institute of Mathematics\\ Polish Academy of Sciences\\ \'Snia\-deckich 8\\ 00-956 Warsaw\\ Poland}
\email{bernard@impan.pl}
\address{Wojciech M. Zajączkowski\\ Institute of Mathematics\\ Polish Academy of Sciences\\ \'Sniadeckich 8\\ 00-956 Warsaw\\ Poland\\ and \\ Institute of Mathematics and Crypto\-logy\\ Military University of Technology\\ Kaliskiego 2\\ 00-908 Warsaw\\ Poland}
\email{wz@impan.pl}
\address{Adam Kubica \\ Faculty of Mathematics and Information Science \\ Warsaw University of Technology\\ Koszykowa 75 \\ Warsaw 00-662\\ Poland \\and 
\\ Institute of Mathematics and Crypto\-logy\\ Military University of Technology\\ Kaliskiego 2\\ 00-908 Warsaw\\ Poland}
\email{a.kubica@mini.pw.edu.pl}

\begin{abstract}
	We present a short and elegant proof of the inequality $\norm{p}_{L_s(\Omega)} \leq c(\Omega) \left(\norm{v}^2_{L_{2s}(\Omega)} + \norm{f}_{L_s(\Omega)}\right)$ for bounded domains $\Omega$ under the slip and Navier boundary conditions. We also show an application of this result for conditional regularity of weak solutions to the Navier-Stokes equations.
\end{abstract}

\maketitle

\section{Introduction}

We consider the initial-boundary value problem for the Navier-Stokes equations
\begin{equation}\label{p1}
	\begin{aligned}
		&v_{,t} + (v\cdot \nabla) v - \nu \triangle v + \nabla p = f & &\text{in $\Omega\times(0,T) =: \Omega^T$},\\
		&\Div v = 0 & &\text{in $\Omega^T$}, \\
		&v\vert_{t = 0} = v(0) & &\text{in $\Omega$},
	\end{aligned}
\end{equation}
either with boundary slip conditions
\begin{equation}\label{p2}
	\begin{aligned}
		\begin{aligned}
			&n \cdot \mathbb{D}(v) \cdot \tau^{\alpha} = 0, \\
			&n \cdot v = 0
		\end{aligned}& &\text{on $\partial \Omega$}
	\end{aligned}
\end{equation}
or with the Navier boundary conditions
\begin{equation}\label{p3}
	\begin{aligned}
		\begin{aligned}
			&\Rot v \times n = 0, \\
			&n \cdot v = 0
		\end{aligned}& &\text{on $\partial \Omega$,}
	\end{aligned}
\end{equation}
where $\Omega \subset \mathbb{R}^3$ is a bounded domain. In case of the boundary slip conditions it is more convenient to write \eqref{p1}$_1$ in the form
\begin{equation}\label{eq52}
	v_{,t} + (v\cdot \nabla)v - \Div \mathbb{T}(v,p) = f.
\end{equation}

To make the above conditions clear let us recall that $n$ and $\tau^{\alpha}$, $\alpha \in \{1,2\}$ are the unit outward normal vector and the unit tangent vectors. By $\mathbb{T}(v,p)$ we mean the stress tensor
\begin{equation*}
	\mathbb{T}(v,p) = \nu\mathbb{D}(v) - p\mathbb{I},
\end{equation*}
where $\mathbb{D}(v)$ is the dilatation tensor, which equals $\frac{1}{2}\left(\nabla v + \nabla^{\perp} v\right)$, $\mathbb{I}$ is the unit matrix and $\nu > 0$ represents the viscosity coefficient.

Note that \eqref{p3} is sometimes referred also as a boundary slip condition whereas \eqref{p2} as the Navier boundary condition. In some cases these conditions coincide (i.e. $\Omega$ is half-space) but in general they differ. In certain cases this difference can be measured in term of the curvature of $\partial \Omega$ (for $\Omega$ of cylindrical type see e.g. \cite[Lemma 6.5]{Nowakowski2012}) but this issue is beyond the scope of our work.

It is well known that for $v(0) \in H^1(\Omega)$ there exists at least one weak solution (see e.g. \cite{Hopf:1950fk}, \cite{Galdi:2000uq}), but the problem of uniqueness and regularity of weak solution in three dimensions remains open.

Our primary interest in \eqref{p1} is an extension of the regularity criterion for the weak solutions onto bounded domains under boundary slip type conditions. One of the basic ideas used in the proof would rely on testing \eqref{p1} with $v\abs{v}^{\theta - 2}$, $\theta \geq 2$. This approach leads to a difficulty related to the estimate for the pressure term. In the whole space or in periodic setting it can be resolved by the application of the Calder\'on-Zygmund theorem (see e.g. \cite{Struwe:2007vn}) to the equation
\begin{equation}\label{eq210}
	-\triangle p = \sum_{i,j=1}^3\frac{\partial^2}{\partial x_i\partial x_j} \left(v_iv_j\right),
\end{equation}
thereby yielding the following estimate
\begin{equation}\label{eq340}
	\norm{p}_{L_s(\Omega)} \leq \norm{v}_{L_{2s}(\Omega)}^2.
\end{equation}
Clearly, in bounded domains \eqref{eq210} must be supplemented with some boundary condition, which at large are difficult or even impossible to establish due to lack of information on $p$ or $\frac{\partial p}{\partial n}$ on $\partial \Omega$ in terms of $v$. One of effective, but restrictive to particular cases remedies that may be exhausted lies in e.g. choosing axially symmetric cylinders with boundary slip conditions (see e.g. \cite[Ch. 3, Lemma 1.1]{Zajaczkowski:2004fk}). One can put another restrictions on the geometry of the domain or on the boundary conditions but at the end we lose certain generality. Therefore, it is reasonable to look for any estimates to for the pressure without the necessity of analyzing \eqref{eq210}.

In this paper we give an alternative proof of the estimate of the form of \eqref{eq340}, which indeed does not rely on \eqref{eq210}. In principle, it is based on an auxiliary Poisson equation with the Neumann boundary conditions. The result reads:

\begin{thm}\label{thm1}
	Suppose that $f \in L_s(\Omega)$ and $v \in L_{2s}(\Omega)$ satisfy (\ref{p1}). Let \eqref{p3} hold and $\Omega$ is bounded and sufficiently regular. Then $p \in L_s(\Omega)$ and
	\begin{equation*}
		\norm{p}_{L_s(\Omega)} \leq c(\Omega)\left(\norm{v}_{L_{2s}(\Omega)}^2 + \norm{f}_{L_s(\Omega)}\right).
	\end{equation*}
	Let now \eqref{p2} hold. If in addition $\nabla v \in L_s(\Omega)$, then
	\begin{equation*}
		\norm{p}_{L_s(\Omega)} \leq c(\Omega)\left(\norm{v}_{L_{2s}(\Omega)}^2 + \norm{\nabla v}_{L_s(\Omega)} + \norm{f}_{L_s(\Omega)}\right).
	\end{equation*}
\end{thm}

\begin{rem}
	The regularity requirement concerning the domain $\Omega$ is related to the Neumann problem for the Poisson equation. For a given $s$ the set $\Omega$ is \emph{sufficiently regular} if for each $g \in L_{s'}(\Omega)$ such that $\int_{\Omega}g\, \ud x = 0$  the following problem
	\begin{equation}\label{np4}
		\begin{aligned}
		&-\triangle \psi = g\,  & &\text{in $\Omega$}, \\
		&n\cdot \nabla \psi = 0 & &\text{on $S$}, \\
		&\int_{\Omega} \psi=0,
	\end{aligned}
	\end{equation}
	has the unique solution $u \in W^2_{s'}(\Omega)$, where $\frac{1}{s}+ \frac{1}{s'}=1$. It holds, for example, if:
	\begin{itemize}
		\item $\partial \Omega \in C^{1,1}$, $s>1$ (see \cite[Lemma 2.4.2.1]{Grisvard:1985vn})
		
		\item $\Omega$ is convex, $s\in [2, \infty)$ (see \cite{Adolfsson:1994uq}),
		
		\item $\Omega\subset \mathbb{R}^{3}$ is a bounded, convex polyhedron and
			\begin{equation*}
				1< \frac{s}{s-1} < \min\left\{3,\frac{2\alpha_{E}}{(2\alpha_{E} - \pi)_{+}} \right\},
			\end{equation*}
		where $\alpha_{E}$ is opening  of the dihedral angle with edge $E$ (see \cite{Mazya:2009fk}),
		\item $\Omega=[0,a]\times[0,b]\times[0,c]$, $s >1$ (the reflection argument),
		
		\item $\Omega = [0,a] \times \Omega'$, where $\Omega' \subset \mathbb{R}^2$ is a bounded set with a smooth boundary, $s > 1$ (the reflection argument).
	\end{itemize}
\end{rem}

Clearly, apart the already mentioned idea, there are different techniques that could be utilized to analyze problem \eqref{p1} with or without relying on \eqref{eq210}. It was a great surprise that they seem to work only in case of the Dirichlet boundary conditions (see \cite{Choe:1998kx}, \cite[Sec. 3]{Berselli:2002ys}, \cite{Zhou:2004fk}, \cite{Kang:2006uq}, \cite{Farwig:2009zr}, \cite{Kim:2010qf}), whereas the boundary slip type conditions were only considered in the half space (see \cite{Bae:2008fk} and \cite{Bae:2008uq}) or in the case of axially symmetric solutions (see \cite{Zajaczkowski:2010lr}). In our work we achieve a little progress. Although the domain we work with is bounded but we assume that it is of cubical type. This kind of restriction, tightly related to the boundary integrals, is removable in many cases (see Remark \ref{rem6}). Our major motivation for investigating the simplest domain follows from intention of keeping the calculations clear and simple. The result reads:

\begin{thm}\label{thm2}
	Let $f\equiv 0$, $T > 0$ and $\Omega := [0,a]\times[0,b]\times[0,c]$ for finite, positive real constants $a$, $b$ and $c$. Suppose that a weak solution $v$ to \eqref{p1} supplemented with either \eqref{p2} or \eqref{p3} satisfies
	\begin{equation*}
		\norm{v}_{L_q(0,T;L_p(\Omega))} < +\infty \qquad \text{where} \qquad \frac{3}{p} + \frac{2}{q} = 1
	\end{equation*}
	for $q < +\infty$.
	Then, $v$ is unique and smooth.
\end{thm}

\begin{rem}
	The assumption $f\equiv 0$ is artificial and can be omitted. It does not change the proof but makes it a little longer.
\end{rem}

\begin{rem}
	For the definition and the proof of existence of weak solutions to \eqref{p1} supplemented with \eqref{p2} or \eqref{p3} see e.g. Introduction in \cite{Zajaczkowski:2005zr}.
\end{rem}

\begin{rem}
	If we drop the assumption on the cubical shape of the domain, the claim of Theorem \ref{thm2} in case $q < + \infty$ is still true. The proof is different, easier but does not base on Theorem \ref{thm1}. We will present its sketch at the end of this work.
\end{rem}

\begin{rem}\label{rem3}
	In cubical domains under both boundary conditions \eqref{p2} and \eqref{p3} the assertion of Theorem \ref{thm1} reads
	\begin{equation*}
		\norm{p}_{L_s(\Omega)} \leq c(\Omega)\left(\norm{v}_{L_{2s}(\Omega)}^2 + \norm{f}_{L_s(\Omega)}\right).
	\end{equation*}	
\end{rem}

Before we move to the next section, let us note that the extension of Serrin condition is mostly studied for the Cauchy problem (see e.g. \cite{Kozono:2004nx}, \cite{Kukavica:2006cr}, \cite{Zhou:2006kx}, \cite{Cao:2008ly}, \cite{Bjorland:2011ve}, \cite{Penel:2011vn}) or for the local-interior regularity (see e.g. \cite{Gustafson:2006dq}), thereby excluding the boundary issues. We do not intend to compare or discuss these improvements. This has been nicely done in several papers. The interested reader we would refer e.g. to \cite{Berselli:2009bh}.

\section{Auxiliary results}

Throughout this article we use the following Young inequality:
\begin{lem}[The Young inequality with small parameter $\kappa$]\label{lem4}
	For any positive $a$ and $b$ the inequality
\begin{equation*}
        ab \leq \kappa a^{\lambda_1} + (\kappa \lambda_1)^{-\frac{\lambda_2}{\lambda_1}}\lambda_2^{-1}b^{\lambda_2}
    \end{equation*}
    holds, where $\kappa > 0$ and
    \begin{equation*}
    	\frac{1}{\lambda_1} + \frac{1}{\lambda_2} = 1, \qquad 1 < \lambda_1, \lambda_2 < + \infty.
    \end{equation*}
\end{lem}

Another useful tool is the imbedding lemma for the space $V^k_2(\Omega^t)$, which is defined as the closure of $\mathcal{C}^{\infty}(\Omega\times(t_0,t_1))$ in the norm
\begin{equation*}
	\norm{u}_{V^k_2(\Omega^t)}^2 = \underset{t\in (t_0,t_1)}{\esssup}\norm{u(t)}_{H^k(\Omega)}^2 \\
	+\left(\int_{t_0}^{t_1}\norm{\nabla u(t)}^2_{H^{k}(\Omega)}\, \ud t\right)^{1/2}.
\end{equation*}
The imbedding lemma reads:
\begin{lem}\label{lem2}
	Suppose that $u \in V_2^0(\Omega^t)$, where $\Omega^t := \Omega\times (t_0,t)$, $t_0 < t \leq t_1$. Then $u \in L_q(t_0,t;L_p(\Omega))$ and
	\begin{equation*}
		\norm{u}_{L_q(t_0,t;L_p(\Omega))} \leq c(p,q,\Omega) \norm{u}_{V^0_2(\Omega^t)}
	\end{equation*}
	holds under the condition $\frac{3}{p} + \frac{2}{q} = \frac{3}{2}$, $2 \leq p \leq 6$.
\end{lem}
Let us emphasize that the constant that appears on the right-hand side does not depend on time.

\begin{lem}[Imbedding theorem]\label{lem10}
	Let $\Omega$ satisfy the cone condition and let $q \geq p$. Set
	\begin{equation*}
		\kappa = 2 - 2r - s - 5\left(\frac{1}{p} - \frac{1}{q}\right) \geq 0.
	\end{equation*}
Then for any function $u \in W^{2,1}_{p}(\Omega^t)$ the inequality
	\begin{equation*}
		\norm{\partial^r_t \Ud_x^s u}_{L_q(\Omega^t)} \leq c_1(p,q,r,s,\Omega) \epsilon^{\kappa}\norm{u}_{W^{2,1}_p(\Omega^t)} + c_2(p,q,r,s,\Omega) \epsilon^{-\kappa + 2s - 2}\norm{u}_{L_p(\Omega^t)}
	\end{equation*}
	holds, where the constants $c_1$ and $c_2$ do not depend on $t$.
\end{lem}
For the proof of the lemma we refer the reader to \cite[Ch.2, \S 3, Lemma 3.3]{lad}.

\section{Proof of Theorem \ref{thm1}}

\begin{proof}
	Let $\psi \in W^{2}_{s'}(\Omega)$ be a unique solution to the following elliptic problem:
	\begin{equation}\label{p4}
		\begin{aligned}
			&-\triangle \psi = p\abs{p}^{s - 2} - \frac{1}{\abs{\Omega}}\int_{\Omega}p\abs{p}^{s - 2}\, \ud x & &\text{in $\Omega$}, \\
			&n\cdot \nabla \psi = 0 & &\text{on $S$}, \\
			&\int_{\Omega}\psi\, \ud x = 0 .
		\end{aligned}
	\end{equation}
	Then the estimate
	\begin{equation}\label{eq2}
		\norm{\psi}_{W^2_{s'}(\Omega)} \leq c(s,\Omega) \norm{p}^{s - 1}_{L_{s}(\Omega)}
	\end{equation}
	holds. Multiplying \eqref{p1}$_1$ by $\nabla \psi$ and integrating over $\Omega$ yields
	\begin{equation}\label{eq1}
		\int_{\Omega} \big(v_{,t} - \nu \triangle v + v\cdot \nabla v + \nabla p\big)\cdot \nabla \psi\, \ud x = \int_{\Omega} f \cdot \nabla \psi\, \ud x.
	\end{equation}
	We have four integral on the left-hand side which need to be estimated. First we see
	\begin{equation}\label{eq50}
		\int_{\Omega} v_{,t} \cdot \nabla \psi \,\ud x = \int_{\Omega}\big(\nabla \cdot (v_{,t}\psi) - \Div v_{,t}\psi\big)\, \ud x = \int_{S} \psi(v_{,t} \cdot n)\,\ud S = 0.
	\end{equation}
	The estimate for the second integral varies in dependence on the boundary conditions. Let us assume \eqref{p3} first. Condition \eqref{p2} will be discussed at the end of the proof. Thus,
	\begin{equation*}
		\int_{\Omega} \triangle v\cdot \nabla \psi\, \ud x = -\int_{\Omega} \Rot\Rot v\cdot \nabla \psi\, \ud x = -\int_{S} \Rot v \times n \cdot \nabla \psi\, \ud S = 0.
	\end{equation*}
	For the third integral we have
	\begin{equation*}
		\int_{\Omega} v_iv_{j,x_i}\psi_{,x_j}\, \ud x = -\int_{\Omega} v_iv_j\psi_{,x_jx_i}\, \ud x + \int_S v_iv_j\psi_{,x_j} n_i\, \ud S \leq \norm{v}^2_{L_{2s}(\Omega)}\norm{\nabla^2\psi}_{L_{s'}(\Omega)},
	\end{equation*}
	where we integrated by parts and utilized  equality (\ref{p3})${}_{2}$. The last term on the left-hand side in \eqref{eq1} is equal to	
	\begin{equation}\label{eq54}
		- \int_{\Omega} p \cdot \triangle \psi\, \ud x + \int_S p \left(n\cdot \nabla \psi\right) \, \ud S = -\int_{\Omega} \abs{p}^{s } \ud x + \frac{1}{\abs{\Omega}} \left(\int_{\Omega} p\abs{p}^{s - 2}\, \ud x\right)\, \int_{\Omega} p\, \ud x = -\norm{p}_{L_s(\Omega)}^s ,
	\end{equation}
	because the boundary integral is equal to zero due to \eqref{p4}$_2$ and $p$ is a distribution determined up to a constant.

	Finally, by the H\"older inequality
	\begin{equation*}
		\int_{\Omega} f\cdot \nabla \psi\, \ud x \leq \norm{f}_{L_s(\Omega)}\norm{\nabla \psi}_{L_{s'}(\Omega)}.
	\end{equation*}
	
	Summing up the above estimates and in view of \eqref{eq2} we obtain	
	\begin{equation*}
		\norm{p}_{L_s(\Omega)}^s \leq \norm{v}^2_{L_{2s}(\Omega)}\norm{\nabla^2\psi}_{L_{s'}(\Omega)} + \norm{f}_{L_s(\Omega)}\norm{\nabla \psi}_{L_{s'}(\Omega)} \leq c(\Omega) \norm{p}^{s - 1}_{L_{s}(\Omega)}\left(\norm{v}_{L_{2s}(\Omega)}^2 + \norm{f}_{L_s(\Omega)}\right).
	\end{equation*}
	Hence
	\begin{equation*}
		\norm{p}_{L_s(\Omega)} \leq c(\Omega)\left(\norm{v}_{L_{2s}(\Omega)}^2 + \norm{f}_{L_s(\Omega)}\right),
	\end{equation*}
	which concludes the proof of the first assertion.
	
	Let now \eqref{p2} hold. Then, instead of \eqref{eq1} we have in light of \eqref{eq52} the identity
	\begin{equation}\label{eq60}
		\int_{\Omega} \big(v_{,t} - \Div \mathbb{T}(v,p) + v\cdot \nabla v \big)\cdot \nabla \psi\, \ud x = \int_{\Omega} f \cdot \nabla \psi\, \ud x.
	\end{equation}
	We need to examine the term involving the Cauchy stress tensor. We see that
	\begin{multline}\label{eq58}
		\int_{\Omega} - \Div \mathbb{T}(v,p) \cdot \nabla \psi\, \ud x \\
		= - \int_S n \cdot \nu\mathbb{D}(v) \cdot \nabla \psi\, \ud S + \int_S p \left(\nabla \psi \cdot n\right)\, \ud S + \int_{\Omega} \nu \mathbb{D}(v) \nabla^2 \psi\, \ud x + \int_{\Omega} p \triangle \psi\, \ud x.
	\end{multline}
	Expressing $\mathbb{D}(v)$ in the basis $n,\tau^{\alpha}$, $\alpha = 1,2$ yields
	\begin{equation*}
		\int_S n \cdot \nu\mathbb{D}(v) \cdot \nabla \psi\, \ud S  = \nu \int_S \left(n \cdot \mathbb{D}(v) \cdot n\right)n \cdot \nabla \psi\, \ud S + \nu \int_S \left(n \cdot \mathbb{D}(v) \cdot \tau^{\alpha}\right)\tau^{\alpha} \cdot \nabla \psi\, \ud S.
	\end{equation*}
	The first integral vanishes due to \eqref{p4}$_2$, whereas the second due to \eqref{p2}. Combining \eqref{eq50}, \eqref{eq54} and \eqref{eq58} we infer from \eqref{eq60} that
	\begin{multline*}
		\norm{p}_{L_s(\Omega)}^s \leq \norm{v}^2_{L_{2s}(\Omega)}\norm{\nabla^2\psi}_{L_{s'}(\Omega)} + \norm{\nabla v}_{L_{s}(\Omega)}\norm{\nabla^2\psi}_{L_{s'}(\Omega)} + \norm{f}_{L_s(\Omega)}\norm{\nabla \psi}_{L_{s'}(\Omega)} \\
			\leq c(\Omega) \norm{p}^{s - 1}_{L_{s}(\Omega)}\left(\norm{v}_{L_{2s}(\Omega)}^2 + \norm{\nabla v}_{L_s(\Omega)} + \norm{f}_{L_s(\Omega)}\right),
	\end{multline*}
	which is our second assertion. The proof is complete.
\end{proof}

\section{Proof of Theorem \ref{thm2}}

\begin{proof}
	We start with multiplying \eqref{p1} by $v\abs{v}^{\theta - 2}$ and integrating over $\Omega$
	\begin{equation}\label{eq3}
	        \frac{1}{\theta}\Dt \int_{\Omega}\abs{v}^{\theta}\, \ud x + \nu\int_\Omega \nabla v \cdot \nabla \left(v\abs{v}^{\theta - 2}\right)\, \ud x = -\int_\Omega \nabla p\cdot v\abs{v}^{\theta - 2}\, \ud x + \int_S  \sum_{i,j = 1}^3v_{j,x_i} \cdot v_j\abs{v}^{\theta - 2} \cdot n_i\, \ud S,
    \end{equation}
    where $\theta > 3$ and the non-linear term vanishes due to
    \begin{multline*}
    	\int_{\Omega} (v\cdot \nabla)v\cdot v\abs{v}^{\theta -2 }\, \ud x = \frac{1}{2}\int_{\Omega} v_i \left(v\cdot v\right)_{,x_i} \abs{v}^{\theta - 2}\, \ud x = \frac{1}{\theta}\int_{\Omega} v_i \left(\abs{v}^2\right)^{\frac{\theta}{2}}_{,x_i}\, \ud x \\
    	= -\frac{1}{\theta}\int_{\Omega} \Div v \abs{v}^{\theta}\, \ud x + \frac{1}{\theta}\int_S \abs{v}^{\theta} \left(v\cdot n\right)\, \ud S = 0.
    \end{multline*}
    Consider first the boundary integral on the right-hand side. On the walls $x_3 = 0$ and $x_3 = c$ the normal vector $n$ equals $(0,0,\mp1)$ and conditions \eqref{p2}, \eqref{p3} imply $v_{1,x_3} = v_{2,x_3} = v_3 = 0$ (see \cite[Lemma 3.1 and its proof]{Zajaczkowski:2005zr} and \cite[Lemma 6.6]{Nowakowski2012}, respectively). Therefore
    \begin{equation*}
    	\int_{S \cap x_3 \in \{0,c\}}  \sum_{i,j = 1}^3v_{j,x_i} \cdot v_j\abs{v}^{\theta - 2} \cdot n_i\, \ud S = 0.
    \end{equation*}
	Following nearly identical reasoning for $x_2 \in \{0,b\}$ and $x_1 \in \{0,a\}$ we conclude that
    \begin{equation*}
    	\int_{S} \sum_{i,j = 1}^3v_{j,x_i} \cdot v_j\abs{v}^{\theta - 2} \cdot n_i\, \ud S = 0.
    \end{equation*}	
    For the second term on the left-hand side in \eqref{eq3} we have
    \begin{equation*}
        \nu\int_\Omega \nabla v \cdot \nabla \left(v\abs{v}^{\theta - 2}\right)\, \ud x = \nu\int_\Omega \abs{\nabla v}^2 \abs{v}^{\theta - 2}\, \ud x + \frac{4\nu(\theta - 2)}{\theta^2}\int_\Omega \abs{\nabla \abs{v}^{\frac{\theta}{2}}}^2 \, \ud x.
    \end{equation*}
    To estimate the term with the pressure we integrate by parts and use (\ref{p1})${}_{2}$ and boundary conditions
    \begin{equation*}
    	-\int_\Omega \nabla p\cdot v\abs{v}^{\theta - 2}\, \ud x = \left(\frac{\theta}{2} - 1\right)\int_{\Omega} p |v|^{\theta-4}v  \cdot \nabla  \abs{v}^{2}\,  \ud x \leq (\theta - 2)\int_{\Omega} \abs{p}\abs{v}^{\frac{\theta}{2} - 1}\abs{\nabla v}\abs{v}^{\frac{\theta}{2} - 1}\, \ud x.
    \end{equation*}
    From the Cauchy inequality we immediately get
    \begin{equation*}
	    (\theta - 2)\int_{\Omega} \abs{p}\abs{v}^{\frac{\theta}{2} - 1}\abs{\nabla v}\abs{v}^{\frac{\theta}{2} - 1}\, \ud x \leq (\theta  - 2) \left(\int_{\Omega} \abs{p}^2 \abs{v}^{\theta - 2}\, \ud x\right)^{\frac{1}{2}} \left(\int_{\Omega} \abs{\nabla v}^2 \abs{v}^{\theta - 2}\,\ud x\right)^{\frac{1}{2}}.
    \end{equation*}

    So far we have obtained
    \begin{multline}\label{eq18}
    	\frac{1}{\theta}\Dt \int_{\Omega}\abs{v}^{\theta}\, \ud x + \frac{4\nu(\theta - 2)}{\theta^2}\int_\Omega \abs{\nabla \abs{v}^{\frac{\theta}{2}}}^2 \, \ud x +  \nu\int_\Omega \abs{\nabla v}^2 \abs{v}^{\theta - 2}\, \ud x \\
    	\leq (\theta  - 2) \left(\int_{\Omega} \abs{p}^2 \abs{v}^{\theta - 2}\, \ud x\right)^{\frac{1}{2}} \left(\int_{\Omega} \abs{\nabla v}^2 \abs{v}^{\theta - 2}\,\ud x\right)^{\frac{1}{2}}.
    \end{multline}
	To estimate the right-hand side we use the H\"older inequality
	\begin{multline*}
		\left(\int_{\Omega} \abs{p}^2 \abs{v}^{\theta - 2}\, \ud x\right)^{\frac{1}{2}} \left(\int_{\Omega} \abs{\nabla v}^2 \abs{v}^{\theta - 2}\,\ud x\right)^{\frac{1}{2}} \\
		\leq \left(\left(\int_{\Omega} \abs{p}^{2\lambda_1}\, \ud x\right)^{\frac{1}{\lambda_1}}\left(\int_{\Omega} \abs{v}^{(\theta - 2)\lambda_2}\, \ud x\right)^{\frac{1}{\lambda_2}}\right)^{\frac{1}{2}} \left(\int_{\Omega} \abs{\nabla v}^2 \abs{v}^{\theta - 2}\,\ud x\right)^{\frac{1}{2}}.
	\end{multline*}
	By Remark \ref{rem3}
	\begin{equation*}
		\left(\int_{\Omega} \abs{p}^{2\lambda_1}\, \ud x\right)^{\frac{1}{2\lambda_1}} \leq c(\Omega) \left(\int_{\Omega}\abs{v}^{4\lambda_1}\, \ud x\right)^{\frac{2}{4\lambda_1}} = c(\Omega)\norm{v}_{L_{4\lambda_1}(\Omega)}^2,
	\end{equation*}
	which combined with \eqref{eq18} yields
    \begin{multline}\label{eq20}
    	\frac{1}{\theta}\Dt \int_{\Omega}\abs{v}^{\theta}\, \ud x + \frac{4\nu(\theta - 2)}{\theta^2}\int_\Omega \abs{\nabla \abs{v}^{\frac{\theta}{2}}}^2 \, \ud x +  \nu\int_\Omega \abs{\nabla v}^2 \abs{v}^{\theta - 2}\, \ud x\\
    	\leq c(\Omega)(\theta - 2) \norm{v}_{L_{4\lambda_1}(\Omega)}^2 \norm{v}_{L_{(\theta -2 )\lambda_2}(\Omega)}^{\frac{\theta - 2}{2}}\left(\int_{\Omega} \abs{\nabla v}^2 \abs{v}^{\theta - 2}\,\ud x\right)^{\frac{1}{2}}.
    \end{multline}
    Due to the imbedding $H^1(\Omega) \hookrightarrow L_6(\Omega)$  and the Poincar\'e inequality (every component of $v$ vanishes on different part of the boundary) we see that
    \begin{equation}\label{eq30}
    	\int_\Omega \abs{\nabla \abs{v}^{\frac{\theta}{2}}}^2\, \ud x \geq c(\Omega) \left(\int_{\Omega} \abs{v}^{\frac{\theta}{2} \cdot 6}\, \ud x\right)^{\frac{1}{6}\cdot 2} = c(\Omega) \norm{v}_{L_{3\theta}(\Omega)}^{\theta}.
    \end{equation}
    Therefore we interpolate $L_{4\lambda_1}(\Omega)$ and $L_{(\theta - 2)\lambda_2}(\Omega)$ between $L_{\theta}(\Omega)$ and $L_{3\theta}(\Omega)$:
    \begin{equation*}
    	\begin{aligned}
			\frac{1}{4\lambda_1} &= \frac{\alpha}{\theta} + \frac{1 - \alpha}{3\theta} = \frac{2\alpha + 1}{3\theta} & &\Leftrightarrow & &\alpha = \frac{1}{2}\left( \frac{3\theta}{4\lambda_1} - 1\right) = \frac{3\theta - 4\lambda_1}{8\lambda_1}, \\
			\frac{1}{(\theta - 2)\lambda_2} &= \frac{\beta}{\theta} + \frac{1 - \beta}{3\theta} = \frac{2\beta + 1}{3\theta} & &\Leftrightarrow & &\beta = \frac{1}{2}\left(\frac{3\theta}{(\theta - 2)\lambda_2} - 1\right) = \frac{3\theta - (\theta - 2)\lambda_2}{2(\theta - 2)\lambda_2}
		\end{aligned}
	\end{equation*}
	and
	\begin{equation*}
		\begin{aligned}
			1 - \alpha &= 1 - \frac{3\theta - 4\lambda_1}{8\lambda_1} = \frac{12\lambda_1 - 3\theta}{8\lambda_1}, \\
			1 - \beta &= 1 - \frac{3\theta - (\theta - 2)\lambda_2}{2(\theta - 2)\lambda_2} = \frac{3(\theta - 2)\lambda_2 - 3\theta}{2(\theta - 2)\lambda_2}.
		\end{aligned}
	\end{equation*}
	Finally
	\begin{equation}\label{eq22}
		\norm{v}_{L_{4\lambda_1}(\Omega)}^2 \norm{v}_{L_{(\theta -2 )\lambda_2}(\Omega)}^{\frac{\theta - 2}{2}} \leq \norm{v}_{L_{\theta}(\Omega)}^{w_1}\norm{v}_{L_{3\theta}(\Omega)}^{w_2},
	\end{equation}
	where
	\begin{multline}\label{eq24}
			w_1 = 2 \cdot \frac{3\theta - 4\lambda_1}{8\lambda_1} + \frac{\theta - 2}{2} \cdot \frac{3\theta - (\theta - 2)\lambda_2}{2(\theta - 2)\lambda_2} = \frac{3\theta}{4\lambda_1} - 1 + \frac{3\theta - (\theta - 2)\lambda_2}{4\lambda_2} = \frac{3\theta}{4} - 1 - \frac{\theta}{4} + \frac{1}{2}\\
			= \frac{\theta}{2} - \frac{1}{2}
	\end{multline}
	and
	\begin{multline}\label{eq26}
		w_2 = 2 \cdot \frac{12\lambda_1 - 3\theta}{8\lambda_1} + \frac{\theta - 2}{2} \cdot \frac{3(\theta - 2)\lambda_2 - 3\theta}{2(\theta - 2)\lambda_2} = 3 - \frac{3\theta}{4\lambda_1} + \frac{3(\theta - 2)\lambda_2 - 3\theta}{4\lambda_2} \\
		= 3 - \frac{3\theta}{4} + \frac{3\theta}{4} - \frac{3}{2} = \frac{3}{2}.
	\end{multline}
	Thus, from \eqref{eq20}, \eqref{eq22}, \eqref{eq24} and \eqref{eq26} it follows
    \begin{multline*}
    	\frac{1}{\theta}\Dt \int_{\Omega}\abs{v}^{\theta}\, \ud x + \frac{4\nu(\theta - 2)}{\theta^2}\int_\Omega \abs{\nabla \abs{v}^{\frac{\theta}{2}}}^2 \, \ud x +  \nu\int_\Omega \abs{\nabla v}^2 \abs{v}^{\theta - 2}\, \ud x  \\
    \leq c(\Omega)(\theta - 2)\norm{v}_{L_{\theta}(\Omega)}^{\frac{\theta}{2} - \frac{1}{2}} \norm{v}_{L_{3\theta}(\Omega)}^{\frac{3}{2}}\left(\int_\Omega \abs{\nabla v}^2 \abs{v}^{\theta - 2}\, \ud x \right)^{\frac{1}{2}}.
    \end{multline*}
    Multiplying by $\theta$ and utilizing \eqref{eq30} in the above inequality gives
    \begin{multline*}
    	\Dt \int_{\Omega}\abs{v}^{\theta}\, \ud x + \frac{4\nu(\theta - 2)}{\theta}\int_\Omega \abs{\nabla \abs{v}^{\frac{\theta}{2}}}^2 \, \ud x + \nu\theta \int_\Omega \abs{\nabla v}^2 \abs{v}^{\theta - 2}\, \ud x  \\
    	\leq c(\Omega)(\theta - 2)\theta\norm{v}_{L_{\theta}(\Omega)}^{\frac{\theta}{2} - \frac{1}{2}} \left(\int_{\Omega} \abs{\nabla \abs{v}^{\frac{\theta}{2}}}^2\,\ud x\right)^{\frac{3}{2\theta}}\left(\int_\Omega \abs{\nabla v}^2 \abs{v}^{\theta - 2}\, \ud x\right)^{\frac{1}{2}} \\
    	\leq  c(\Omega)(\theta - 2)\theta^{2} 2^{-\frac{3}{\theta}} \norm{v}_{L_{\theta}(\Omega)}^{\frac{\theta}{2} - \frac{1}{2}} \left(\int_\Omega \abs{\nabla v}^2 \abs{v}^{\theta - 2}\, \ud x\right)^{\frac{3}{2\theta} + \frac{1}{2}}.
    \end{multline*}

    Utilizing the Young inequality (see Lemma \ref{lem4})  we obtain
	\begin{multline}\label{eq28}
    	\Dt \int_{\Omega}\abs{v}^{\theta}\, \ud x + \frac{4\nu(\theta - 2)}{\theta}\int_\Omega \abs{\nabla \abs{v}^{\frac{\theta}{2}}}^2 \, \ud x + \nu\theta \int_\Omega \abs{\nabla v}^2 \abs{v}^{\theta - 2}\, \ud x \\
    	\leq \kappa \left( \int_\Omega \abs{\nabla v}^2 \abs{v}^{\theta - 2}\, \ud x \right)^{\left(\frac{3}{2\theta} + \frac{1}{2}\right)\gamma_1} + \left(\frac{1}{\kappa\gamma_1}\right)^{\frac{\gamma_2}{\gamma_1}} \frac{1}{\gamma_2} \left( c(\Omega)(\theta - 2)\theta^{2} 2^{-\frac{3}{\theta}})\right)^{\gamma_2}\norm{v}_{L_{\theta}(\Omega)}^{\left(\frac{\theta}{2} - \frac{1}{2}\right)\gamma_2}.
    \end{multline}
    Now we chose $\gamma_1$ so it satisfies
    \begin{equation*}
    	\left(\frac{3}{2\theta} + \frac{1}{2}\right)\gamma_1 = 1 \qquad \Leftrightarrow \qquad \gamma_1 = \frac{2}{\frac{3}{\theta} + 1} = \frac{2\theta}{3 + \theta}.
    \end{equation*}
    Thus
    \begin{equation*}
    	\gamma_2 = \frac{\gamma_1}{\gamma_1 - 1} = \frac{\frac{2\theta}{3 + \theta}}{\frac{2\theta}{3 + \theta} - 1} = \frac{2\theta}{3 + \theta}\cdot \frac{3 + \theta}{2\theta - 3 - \theta} = \frac{2\theta}{\theta - 3}.
    \end{equation*}
    Hence
    \begin{equation*}
		\left(\frac{\theta}{2} - \frac{1}{2}\right)\gamma_2 = \left(\frac{\theta}{2} - \frac{1}{2}\right)\cdot \frac{2\theta}{3 + \theta} = \frac{\theta(\theta - 1)}{\theta - 3}.
    \end{equation*}
    and
   	\begin{equation}\label{eq44}
    	\Dt \int_{\Omega}\abs{v}^{\theta}\, \ud x + \nu\int_\Omega \abs{\nabla \abs{v}^{\frac{\theta}{2}}}^2 \, \ud x + \nu \int_\Omega \abs{\nabla v}^2 \abs{v}^{\theta - 2}\, \ud x\leq c(\nu,\theta,\Omega)\norm{v}_{L_{\theta}(\Omega)}^{\frac{\theta(\theta - 1)}{\theta - 3}}.
    \end{equation}
    Since
    \begin{equation*}
    	\frac{\theta(\theta - 1)}{\theta - 3} = \theta \left(1 + \frac{2}{\theta - 3}\right)
    \end{equation*}
    we put $\theta=p$ and using the assumption on $v$ (then $q$ is equal to $\frac{2p}{p - 3}$) we may apply the Gronwall inequality
    \begin{equation*}
    	\sup_{0 \leq t \leq T} \norm{v(t)}_{L_{p}(\Omega)}^{p} \leq \exp\left(c(\nu,p,\Omega)\norm{v}_{L_{q}(0,T;L_{p}(\Omega))}^{q}\right) \norm{v(0)}^{p}_{L_{p}(\Omega)}.
    \end{equation*}
    Integrating \eqref{eq44} with respect to $t$ gives
    \begin{multline*}
    	\sup_{0\leq t \leq T} \norm{v(t)}^{p}_{L_{p(\Omega)}} + \nu\int_{\Omega^T} \abs{\nabla \abs{v}^{\frac{p}{2}}}^2 \, \ud x\, \ud t\leq c(\nu,p,\Omega)\int_0^T\norm{v(t)}_{L_{p}(\Omega)}^{p+q}\, \ud t + \norm{v(0)}^{p}_{L_{p}(\Omega)} \\
    	\leq c(\nu,p,\Omega)  \exp\left(c(\nu,p,\Omega)\norm{v}_{L_{q}(0,T;L_{p}(\Omega))}^{q}\right)\norm{v}_{L_{q}(0,T;L_{p}(\Omega))}^{q} \norm{v(0)}_{L_{p}(\Omega)}^p + \norm{v(0)}^{p}_{L_{p}(\Omega)}.
    \end{multline*}
     By Lemma \ref{lem2} we get
	\begin{multline}\label{eq130}
		\norm{v}_{L_{\frac{5}{3}p}(\Omega^t)} \leq \left[c(\nu,p,\Omega) \exp\left(c(\nu,p,\Omega)\norm{v}_{L_{q}(0,T;L_{p}(\Omega))}^{q}\right) \norm{v}_{L_{q}(0,T;L_{p}(\Omega))}^{q} + 1\right]\norm{v(0)}_{L_{p}(\Omega)} \\
		=: c(\text{data}).
	\end{multline}
	In view of the classical theory (see e.g. \cite{Solonnikov:1964uq}, \cite{Solonnikov:1976kx}, \cite{Solonnikov:1977vn}, \cite{Solonnikov:1990fk} and recently \cite{wz80}) we infer (see Remark \ref{rem4})
	\begin{equation*}
		\norm{v}_{W^{2,1}_s(\Omega^t)} + \norm{\nabla p}_{L_s(\Omega^t)} \leq c(s,\nu,\Omega)\norm{(v \cdot \nabla)v}_{L_s(\Omega^t)}  + \norm{v(0)}_{W^{2 - \frac{2}{s}}_s(\Omega)}.
	\end{equation*}
	for $t \in (0,T)$. By the H\"older inequality
	\begin{equation*}
	\norm{(v \cdot \nabla)v}_{L_s(\Omega^t)} \leq \norm{v}_{L_{\frac{5}{3}p}(\Omega^t)} \norm{\nabla v}_{L_r(\Omega^t)},
	\end{equation*}
	where
	\begin{equation*}
		\frac{1}{\frac{5}{3}p} + \frac{1}{r} = \frac{1}{s}.
	\end{equation*}
	Lemma \ref{lem10} yields
	\begin{equation*}
		\norm{\nabla v}_{L_r(\Omega^t)} \leq \epsilon^{\kappa} \norm{v}_{W^{2,1}_s(\Omega^t)} + \epsilon^{-\kappa} \norm{v}_{L_s(\Omega^t)},
	\end{equation*}
	where
	\begin{equation*}
		\kappa = 2 - 1 - 5\left(\frac{1}{s} - \frac{1}{r}\right) = 1 - \frac{3}{p} > 0 \qquad \Leftrightarrow \qquad p > 3,
	\end{equation*}
	thereby
	\begin{equation*}
		\norm{v}_{W^{2,1}_s(\Omega^t)} + \norm{\nabla p}_{L_s(\Omega^t)} \leq c(\text{data}) \norm{v}_{L_s(\Omega^t)} +  \norm{v(0)}_{W^{2 - \frac{2}{s}}_s(\Omega)}.
	\end{equation*}
	For $s=p$ the right hand side is finite, thus  $v$ and $p$ are  smooth provided  $v(0)$ is smooth. This completes the proof.
\end{proof}

	\begin{center}
		\includegraphics[scale=0.9]{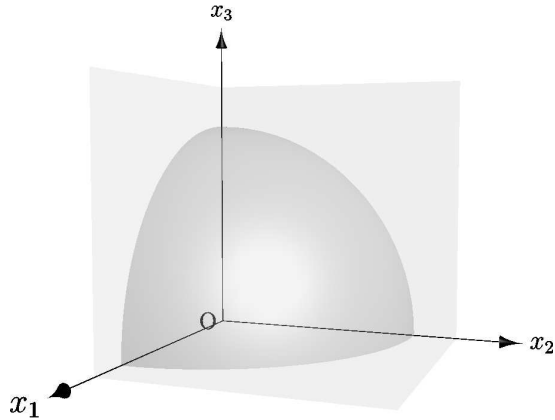}
		\captionof{figure}{The localized problem near the corner $O = (0,0,0)$.}\label{fig:2}
	\end{center}
	
\begin{rem}\label{rem4}
	At the end of the proof of Theorem \ref{thm2} we used some references to classical theory concerning the regularity of the Stokes system under boundary slip conditions. One of the assumptions in these results is certain smoothness of the boundary (roughly speaking: the higher regularity the higher boundary smoothness). In our case we deal with domains of cubical type, which have corners. Nevertheless, the classical theory holds because we can localize the problem near corners and due to either \eqref{p2} or \eqref{p3} reflect it outside the cube. For example, let us consider the corner at $O = (0,0,0)$ (see Figure \ref{fig:2}).

	 As we saw in the beginning of the proof of Theorem \ref{thm2} we have on the wall $x_3 = 0$ the equality $v_{1,x_3} = v_{2,x_3} = v_3 = 0$, which suggests the reflection
	\begin{equation*}
		\check v (x) = \begin{cases}  \bar{v}(x) & x_3 \in \overline{\supp\zeta\cap \Omega}, \\
													( \bar{v}'(\bar x),-\bar{v}_3(\bar x)) & x_3 \leq 0,
						\end{cases}
	\end{equation*}
	where $\bar x = (x',-x_3)$ (see Figure \ref{fig:3}).
	\begin{figure}	
		\centering
		\subfloat[Reflection with respect to $x_3 = 0$.]{
			\label{fig:3}
		 	\includegraphics[scale=0.9]{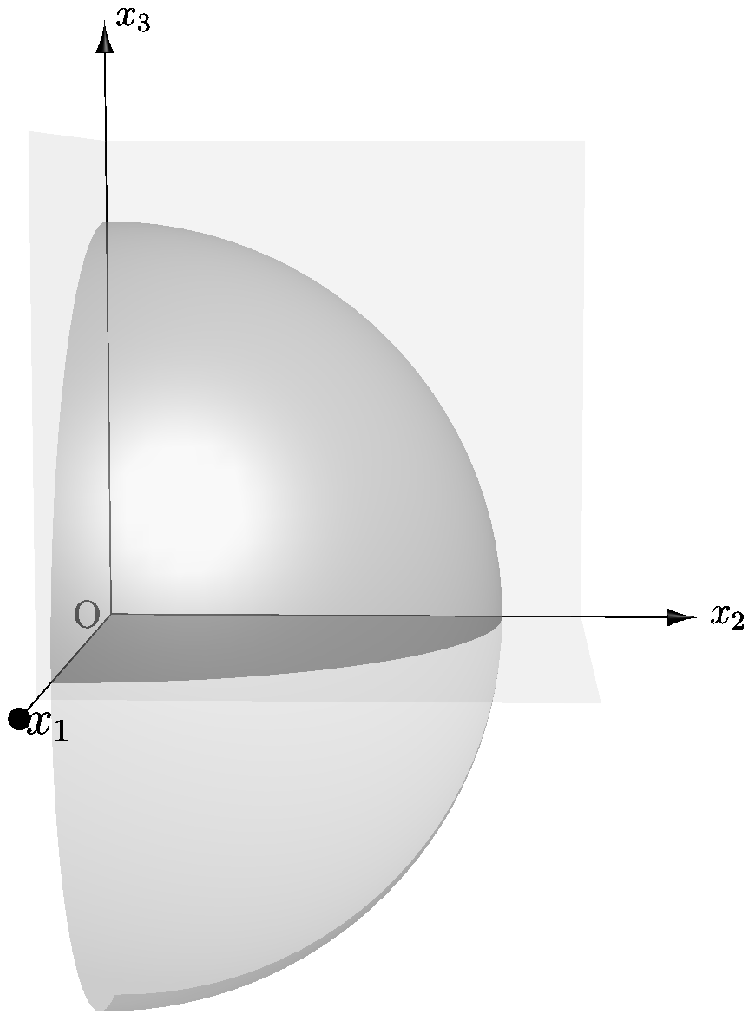}} %
		\subfloat[Reflectiom with respect to $x_3 = 0$ and $x_1 = 0$.]{
			\label{fig:4}
			\includegraphics{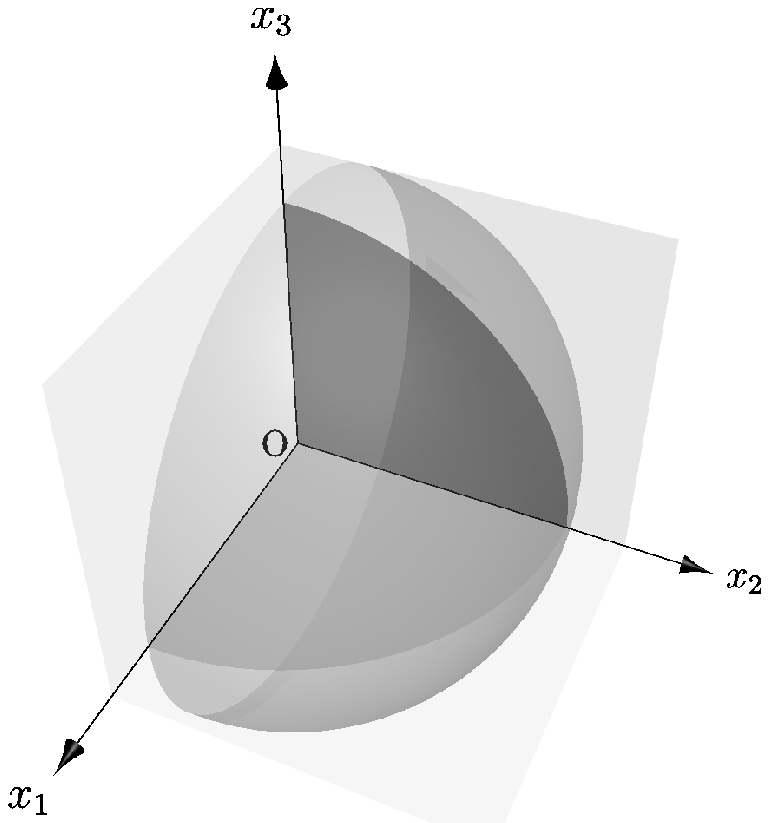}}
		\caption{Reflections of the localized problem.}
		\label{fig:reflections}		
	\end{figure}
	By $\supp \zeta$ we denote the support of the cut-off function $\zeta$, and $\bar v$ denotes $v$ localized to $\supp \zeta$, i.e. $\bar v = v \zeta$. Similarly, since $f = 0$ we immediately get that $\frac{\partial p}{\partial n} = 0$ on each part of the boundary. This implies that the reflection with respect to $x_3$ preserves the Stokes system. Now, to get the problem in the half-space we need one more reflection (see Figure \ref{fig:4}). Observe that on $x_1 = 0$ we have $v_1 = v_{2,x_1} = v_{3,x_1} = 0$ and $\frac{\partial p}{\partial n} = 0$, so we introduce
	\begin{equation*}
		\check{\check{ v}} (x) = \begin{cases} \check v(x) & x_1 \in \overline{\supp \check v}, \\
									(- \check v_1(\bar{\bar x}),\check v_2(\bar{\bar x}),\check v_3(\bar{\bar x})) & x_1 \leq 0,
								\end{cases}
	\end{equation*}
	where $\bar{\bar x} = (-x_1,x_2,x_3)$. Now we see that $\check{\check{v}}$ is defined in the half-space $x_2 \geq 0$ and the Stokes system is preserved.
\end{rem}

\begin{rem}\label{rem6}
	We have already mentioned that the assumption on the cubical shape of the domain can be relaxed. This motivation follows from \eqref{eq3}, where the appearing boundary integral can be written in the form
	\begin{equation}\label{eq310}
		\int_S  \sum_{i,j = 1}^3v_{j,x_i} \cdot v_j\abs{v}^{\theta - 2} \cdot n_i\, \ud S = \int_S \abs{v}^{\theta - 2} \left(\Rot v \times n\right)\cdot v\, \ud S - \int_S \sum_{i,j = 1}^3 v_{i,x_j} n_iv_j\, \ud S.
	\end{equation}
	We see that under \eqref{p2} the first integral on the right-hand side vanishes.
	
	\begin{figure}
		\centering
		\includegraphics{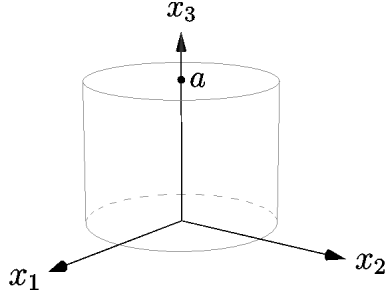}
		\captionof{figure}{The set $\Omega$.}\label{fig:1}		
	\end{figure}

	To eliminate the second integral we impose that $\Omega$ is of cylindrical type, parallel to the $x_3$ axis with convex cross section (see Figure \ref{fig:1}). Denoting the side boundary by $S_1$, the bottom and the top of the cylinder (perpendicular to $x_3$) by $S_2$, the normal unit vector and the tangent unit vectors by $n$, $\tau^{\alpha}$, $\alpha = 1,2$, respectively, we easily establish (see e.g. Introduction in \cite{Zajaczkowski:2005fk}) that
	\begin{equation}\label{p44}
		\begin{aligned}
			&n\vert_{S_1} = \frac{1}{\abs{\nabla \varphi}}(\varphi_{,x_1},\varphi_{,x_2},0) & &\tau^1\vert_{S_1} = \frac{1}{\abs{\nabla \varphi}}(-\varphi_{,x_2},\varphi_{,x_1},0) & & \tau^2\vert_{S_1} = (0,0,1)\\
			&n\vert_{S_2} = \left(0,0,\frac{a}{\abs{a}}\right) & &\tau^1\vert_{S_2} = (1,0,0) & &\tau^2\vert_{S_2} = (0,1,0),
		\end{aligned}
	\end{equation}
	where $\varphi(x_1,x_2) = c_0$ is a sufficiently smooth, convex, closed curve in the plane $x_3 = \operatorname{const}$.
		
	Since $n$ does not depend on $x_3$ on $S_1$ we get in view of \eqref{p44} that
	\begin{equation*}
		\int_{S_1} \sum_{i,j = 1}^3 v_{i,x_j} n_iv_j\, \ud S = \frac{\left(v\cdot \tau_1\right)^2}{\abs{\nabla \varphi}^3} \left(\tau^1_1n_{1,x_1}\tau^1_1 + \tau^1_1n_{1,x_2}\tau^1_2 + \tau^1_2n_{2,x_1}\tau^1_1 + \tau^1_2n_{2,x_2}\tau^1_2\right) = (v \cdot \tau_1)^2 \cdot \kappa,
	\end{equation*}
	where $\kappa$ is the curvature of $\varphi$. On $S_2$ we immediately see that $v_3 = v_{3,x_1} = v_{3,x_2} = 0$. Thus, \eqref{eq310} is negative and can be safely removed from \eqref{eq3}.
	
	For further geometrical considerations of the last term on the right-hand side in \eqref{eq310} we would refer the reader to e.g. \cite[Section 2]{Watanabe:2003fk}.
	
\end{rem}

\begin{rem}
	If we do not use the estimate for the pressure from Theorem \ref{thm1}, then we proceed as follows. 
	First, we multiply \eqref{p1}$_{1,2}$, \eqref{p2} and \eqref{p3} by $\eta_k(t)$, $k \in \mathbb{N}$, where
	\begin{equation*}
		\eta_k(t) =	\begin{cases}
								1 & \text{for } t \in \big(kT,(k + 1)T\big), \\
								0 & \text{for } t \leq (k - 1)T
							\end{cases}		
	\end{equation*}
	with the properties $\eta_k \in \mathcal{C}^{\infty}_c(0,\infty)$ and $\frac{\ud }{\ud t} \eta_k(t) \leq \frac{1}{T}$. Denoting $\bar v = v \eta_k$ (we omit $k$ for clarity) we see that \eqref{p1} becomes	
	\begin{equation}\label{p11}
		\begin{aligned}
			&\bar v_{,t} + (v\cdot \nabla) \bar v - \nu \triangle \bar v + \nabla \bar p = \bar f - v \eta_{,t} =: \bar F& &\text{in $\Omega\times\big((k - 1)T,(k + 1)T\big) =: \Omega^{kT}$},\\
			&\Div \bar v = 0 & &\text{in $\Omega^{kT}$}, \\
			&v\vert_{t = (k - 1)T} = 0 & &\text{in $\Omega$}.
		\end{aligned}
	\end{equation}
	and for \eqref{p2} and \eqref{p3} we have
	\begin{equation}\label{p22}
		\begin{aligned}
			\begin{aligned}
				&n \cdot \mathbb{D}(\bar v) \cdot \tau_{\alpha} = 0, \\
				&n \cdot \bar v = 0
			\end{aligned}& &\text{on $\partial \Omega$}
		\end{aligned}
	\end{equation}
	and
	\begin{equation}\label{p33}
		\begin{aligned}
			\begin{aligned}
				&\Rot \bar v \times n = 0, \\
				&n \cdot \bar v = 0
			\end{aligned}& &\text{on $\partial \Omega$.}
		\end{aligned}
	\end{equation}
	By similar reasoning as in \cite{Solonnikov:2002fk} we get
	\begin{multline*}
		\norm{\bar v}_{W^{2,1}_{p,q}(\Omega^{kT})} + \norm{\nabla p}_{L_q((k-1)T,(k + 1)T;L_p(\Omega))} \\
		\leq c(\nu,p,q,T,\Omega)\left(\norm{(v\cdot \nabla)\bar v}_{L_q((k-1)T,(k + 1)T;L_p(\Omega))} + \norm{\bar F}_{L_q((k-1)T,(k + 1)T;L_p(\Omega))}\right).
	\end{multline*}
	The H\"older inequality implies that
	\begin{equation*}
		\norm{(v\cdot \nabla)\bar v}_{L_q((k-1)T,(k + 1)T;L_p(\Omega))} \leq \norm{v}_{L_s((k-1)T,(k + 1)T;L_r(\Omega))}\norm{\nabla \bar v}_{L_{\beta}((k-1)T,(k + 1)T;L_{\alpha}(\Omega))},
	\end{equation*}
	where
	\begin{equation}\label{eq190}
		\frac{1}{r} + \frac{1}{\alpha} = \frac{1}{p} \qquad \text{and} \qquad \frac{1}{s} + \frac{1}{\beta} = \frac{1}{q}.
	\end{equation}
	The imbedding $W^{2,1}_{p,q}(\Omega^{kT}) \hookrightarrow {L_q((k-1)T,(k + 1)T;L_p(\Omega))}$ holds (see e.g. \cite[Ch.3, \S 10.2]{bes}) provided
	\begin{equation*}
		\left(\frac{1}{p} - \frac{1}{\alpha}\right)\frac{3}{2} + \frac{1}{2} + \left(\frac{1}{q} - \frac{1}{\beta}\right) = 1,
	\end{equation*}
	which in view of \eqref{eq190} is equivalent to
	\begin{equation*}
		\frac{3}{r} + \frac{2}{s} = 1.
	\end{equation*}
	Thus,
	\begin{multline}\label{eq200}
		\norm{\bar v}_{W^{2,1}_{p,q}(\Omega^{kT})} + \norm{\nabla p}_{L_q((k-1)T,(k + 1)T;L_p(\Omega))} \\
		\leq c(\nu,p,q,T,\Omega) \left(\norm{v}_{L_s((k-1)T,(k + 1)T;L_r(\Omega))} +\frac{c(initial\ data)}{T} + \norm{\bar f}_{L_q((k-1)T,(k + 1)T;L_p(\Omega))}\right)
	\end{multline}
	for $T$ small enough.
	
	By the same interpolation argument we deduce that the solution to \eqref{p1} satisfies
	\begin{multline*}
		\norm{v}_{W^{2,1}_{p,q}(\Omega^{T})} + \norm{\nabla p}_{L_q(0,T;L_p(\Omega))} \\
		\leq c(\nu,p,q,T,\Omega) \left(\norm{v}_{L_s(0,T;L_r(\Omega))} + \norm{f}_{L_q(0,T;L_p(\Omega))}\right) + \norm{v(0)}_{W^{2 - \frac{2}{p}}(\Omega)}
	\end{multline*}
	for $T$ small enough. Thus, combining the above inequality with \eqref{eq200} and summing over $k$ yields
	\begin{equation*}
		\norm{v}_{W^{2,1}_{p,q}(\Omega^{T})} + \norm{\nabla p}_{L_q(0,T;L_p(\Omega))}  \leq c(initial\ and\ external\ data)
	\end{equation*}
	for arbitrary large $T < +\infty$. Now, the classical theory yields smoothness of $v$ and $p$.
\end{rem}

\bibliographystyle{amsalpha}
\bibliography{bibliography}

\providecommand{\MR}[1]{}
\providecommand{\bysame}{\leavevmode\hbox to3em{\hrulefill}\thinspace}
\providecommand{\MR}{\relax\ifhmode\unskip\space\fi MR }
\providecommand{\MRhref}[2]{%
  \href{http://www.ams.org/mathscinet-getitem?mr=#1}{#2}
}
\providecommand{\href}[2]{#2}
\begin{thebibliography}{GKT06}

\bibitem[AJ94]{Adolfsson:1994uq}
V.~Adolfsson and D.~Jerison, \emph{{$L^p$}-integrability of the second order
  derivatives for the {N}eumann problem in convex domains}, Indiana Univ. Math.
  J. \textbf{43} (1994), no.~4, 1123--1138.

\bibitem[BCJ08]{Bae:2008uq}
H.-O. Bae, H.J. Choe, and B.J. Jin, \emph{Pressure representation and boundary
  regularity of the {N}avier-{S}tokes equations with slip boundary condition},
  J. Differential Equations \textbf{244} (2008), no.~11, 2741--2763.

\bibitem[Ber09]{Berselli:2009bh}
L.C. Berselli, \emph{Some criteria concerning the vorticity and the problem of
  global regularity for the 3{D} {N}avier-{S}tokes equations}, Ann. Univ.
  Ferrara Sez. VII Sci. Mat. \textbf{55} (2009), no.~2, 209--224.

\bibitem[BG02]{Berselli:2002ys}
L.C. Berselli and G.P. Galdi, \emph{Regularity criteria involving the pressure
  for the weak solutions to the {N}avier-{S}tokes equations}, Proc. Amer. Math.
  Soc. \textbf{130} (2002), no.~12, 3585--3595.

\bibitem[BIN78]{bes}
O.V. Besov, V.P. Il'in, and S.M. Nikol'ski{\u \i}, \emph{Integral
  representations of functions and imbedding theorems. {V}ol. {I}}, V. H.
  Winston \& Sons, Washington, D.C., 1978, Translated from the Russian, Scripta
  Series in Mathematics, Edited by Mitchell H. Taibleson.

\bibitem[BJ08]{Bae:2008fk}
H.-O. Bae and B.J. Jin, \emph{Regularity for the {N}avier-{S}tokes equations
  with slip boundary condition}, Proc. Amer. Math. Soc. \textbf{136} (2008),
  no.~7, 2439--2443.

\bibitem[BV11]{Bjorland:2011ve}
C.~Bjorland and A.~Vasseur, \emph{Weak in space, log in time improvement of the
  {L}ady\v zenskaja-{P}rodi-{S}errin criteria}, J. Math. Fluid Mech.
  \textbf{13} (2011), no.~2, 259--269.

\bibitem[Cho98]{Choe:1998kx}
H.J. Choe, \emph{Boundary regularity of weak solutions of the {N}avier-{S}tokes
  equations}, J. Differential Equations \textbf{149} (1998), no.~2, 211--247.

\bibitem[CT08]{Cao:2008ly}
C.~Cao and E.S. Titi, \emph{Regularity criteria for the three-dimensional
  {N}avier-{S}tokes equations}, Indiana Univ. Math. J. \textbf{57} (2008),
  no.~6, 2643--2661.

\bibitem[FKS09]{Farwig:2009zr}
R.~Farwig, H.~Kozono, and H.~Sohr, \emph{Energy-based regularity criteria for
  the {N}avier-{S}tokes equations}, J. Math. Fluid Mech. \textbf{11} (2009),
  no.~3, 428--442.

\bibitem[Gal00]{Galdi:2000uq}
G.P. Galdi, \emph{An introduction to the {N}avier-{S}tokes initial-boundary
  value problem}, Fundamental directions in mathematical fluid mechanics, Adv.
  Math. Fluid Mech., Birkh{\"a}user, Basel, 2000, pp.~1--70.

\bibitem[GKT06]{Gustafson:2006dq}
S.~Gustafson, K.~Kang, and T.-P. Tsai, \emph{Regularity criteria for suitable
  weak solutions of the {N}avier-{S}tokes equations near the boundary}, J.
  Differential Equations \textbf{226} (2006), no.~2, 594--618.

\bibitem[Gri85]{Grisvard:1985vn}
P.~Grisvard, \emph{Elliptic problems in nonsmooth domains}, Monographs and
  Studies in Mathematics, vol.~24, Pitman (Advanced Publishing Program),
  Boston, MA, 1985.

\bibitem[Hop50]{Hopf:1950fk}
E.~Hopf, \emph{{{\"U}ber die Anfangswertaufgabe f{\"u}r die hydrodynamischen
  Grundgleichungen.}}, Math. Nachr. \textbf{4} (1950), no.~1, 213--231.

\bibitem[Kim10]{Kim:2010qf}
J.-M. Kim, \emph{On regularity criteria of the {N}avier-{S}tokes equations in
  bounded domains}, J. Math. Phys. \textbf{51} (2010), no.~5, 053102, 7.

\bibitem[KL06]{Kang:2006uq}
Kyungkeun Kang and Jihoon Lee, \emph{On regularity criteria in conjunction with
  the pressure of {N}avier-{S}tokes equations}, Int. Math. Res. Not. (2006),
  Art. ID 80762, 25. \MR{2250016 (2007j:35159)}

\bibitem[KS04]{Kozono:2004nx}
H.~Kozono and Y.~Shimada, \emph{Bilinear estimates in homogeneous
  {T}riebel-{L}izorkin spaces and the {N}avier-{S}tokes equations}, Math.
  Nachr. \textbf{276} (2004), 63--74.

\bibitem[KZ06]{Kukavica:2006cr}
I.~Kukavica and M.~Ziane, \emph{One component regularity for the
  {N}avier-{S}tokes equations}, Nonlinearity \textbf{19} (2006), no.~2,
  453--469.

\bibitem[LSU67]{lad}
O.A. Lady{\v z}enskaja, V.A. Solonnikov, and N.N. Ural'ceva, \emph{Linear and
  quasilinear equations of parabolic type}, Translated from the Russian by S.
  Smith. Translations of Mathematical Monographs, Vol. 23, American
  Mathematical Society, Providence, R.I., 1967.

\bibitem[Maz09]{Mazya:2009fk}
V.~Maz'ya, \emph{On the boundedness of first derivatives for solutions to the
  {N}eumann-{L}aplace problem in a convex domain}, J. Math. Sci. (N. Y.)
  \textbf{159} (2009), no.~1, 104--112, Problems in mathematical analysis. No.
  40.

\bibitem[Now13]{Nowakowski2012}
B.~Nowakowski, \emph{Large time existence of strong solutions to micropolar
  equations in cylindrical domains}, Nonlinear Anal. Real World Appl.
  \textbf{14} (2013), no.~1, 635--660.

\bibitem[PP11]{Penel:2011vn}
P.~Penel and M.~Pokorn{{\'y}}, \emph{On anisotropic regularity criteria for the
  solutions to 3{D} {N}avier-{S}tokes equations}, J. Math. Fluid Mech.
  \textbf{13} (2011), no.~3, 341--353.

\bibitem[Sol64]{Solonnikov:1964uq}
V.~A. Solonnikov, \emph{Estimates for solutions of a non-stationary linearized
  system of {N}avier-{S}tokes equations}, Trudy Mat. Inst. Steklov. \textbf{70}
  (1964), 213--317.

\bibitem[Sol76]{Solonnikov:1976kx}
\bysame, \emph{Estimates of the solution of a certain initial-boundary value
  problem for a linear nonstationary system of {N}avier-{S}tokes equations},
  Zap. Nau\v cn. Sem. Leningrad. Otdel Mat. Inst. Steklov. (LOMI) \textbf{59}
  (1976), 178--254, 257, Boundary value problems of mathematical physics and
  related questions in the theory of functions, 9.

\bibitem[Sol77]{Solonnikov:1977vn}
\bysame, \emph{The solvability of the second initial-boundary value problem for
  a linear nonstationary system of {N}avier-{S}tokes equations}, Zap. Nau\v cn.
  Sem. Leningrad. Otdel. Mat. Inst. Steklov. (LOMI) \textbf{69} (1977),
  200--218, 277, Boundary value problems of mathematical physics and related
  questions in the theory of functions, 10.

\bibitem[Sol90]{Solonnikov:1990fk}
\bysame, \emph{An initial-boundary value problem for a {S}tokes system that
  arises in the study of a problem with a free boundary}, Trudy Mat. Inst.
  Steklov. \textbf{188} (1990), 150--188, 192, Translated in Proc. Steklov
  Inst. Math. {{\bf{1}}991}, no. 3, 191--239, Boundary value problems of
  mathematical physics, 14 (Russian).

\bibitem[Sol02]{Solonnikov:2002fk}
\bysame, \emph{Estimates of solutions of the {S}tokes equations in {S}. {L}.\
  {S}obolev spaces with a mixed norm}, Zap. Nauchn. Sem. S.-Peterburg. Otdel.
  Mat. Inst. Steklov. (POMI) \textbf{288} (2002), no.~Kraev. Zadachi Mat. Fiz.
  i Smezh. Vopr. Teor. Funkts. 32, 204--231, 273--274.

\bibitem[Str07]{Struwe:2007vn}
M.~Struwe, \emph{On a {S}errin-type regularity criterion for the
  {N}avier-{S}tokes equations in terms of the pressure}, J. Math. Fluid Mech.
  \textbf{9} (2007), no.~2, 235--242.

\bibitem[Wat03]{Watanabe:2003fk}
J.~Watanabe, \emph{On incompressible viscous fluid flows with slip boundary
  conditions}, Proceedings of the 6th {J}apan-{C}hina {J}oint {S}eminar on
  {N}umerical {M}athematics ({T}sukuba, 2002), vol. 159, 2003, pp.~161--172.

\bibitem[Zaj04]{Zajaczkowski:2004fk}
W.M. Zaj{\k a}czkowski, \emph{Global special regular solutions to the
  {N}avier-{S}tokes equations in a cylindrical domain under boundary slip
  conditions}, {GAKUTO International Series, Mathematical Sciences and
  Applications}, no.~21, Gakk{\=o}tosho, 2004.

\bibitem[Zaj05a]{Zajaczkowski:2005fk}
\bysame, \emph{Global regular nonstationary flow for the {N}avier-{S}tokes
  equations in a cylindrical pipe}, Topol. Methods Nonlinear Anal. \textbf{26}
  (2005), no.~2, 221--286.

\bibitem[Zaj05b]{Zajaczkowski:2005zr}
\bysame, \emph{Long time existence of regular solutions to navier-stokes
  equations in cylindrical domains under boundary slip conditions.}, Stud.
  Math. \textbf{169} (2005), no.~3, 243--285.

\bibitem[Zaj10]{Zajaczkowski:2010lr}
\bysame, \emph{A regularity criterion for axially symmetric solutions to the
  {N}avier-{S}tokes equations}, Zap. Nauchn. Sem. S.-Peterburg. Otdel. Mat.
  Inst. Steklov. (POMI) \textbf{385} (2010), no.~Kraevye Zadachi
  Matematicheskoi Fiziki i Smezhnye Voprosy Teorii Funktsii. 41, 54--68, 234.

\bibitem[Zaj11]{wz80}
\bysame, \emph{Nonstationary stokes system in {S}obolev-{S}lobodetski spaces},
  Submitted.

\bibitem[Zho04]{Zhou:2004fk}
Y.~Zhou, \emph{Regularity criteria in terms of pressure for the 3-{D}
  {N}avier-{S}tokes equations in a generic domain}, Math. Ann. \textbf{328}
  (2004), no.~1-2, 173--192.

\bibitem[Zho06]{Zhou:2006kx}
\bysame, \emph{On regularity criteria in terms of pressure for the
  {N}avier-{S}tokes equations in {${\mathbb R}^3$}}, Proc. Amer. Math. Soc.
  \textbf{134} (2006), no.~1, 149--156 (electronic).

\end{thebibliography}

\end{document}